\documentclass[10pt,twoside]{article}
\usepackage{amssymb,amsmath,amsthm,enumerate,geometry}
\input amssym.def
\newtheorem{theorem}{Theorem}[section]
\newtheorem{lemma}{Lemma}[section]
\newtheorem{proposition}{Proposition}[section]
\newtheorem{remark}{Remark}[section]
\newtheorem{corollary}{Corollary}[section]
\newtheorem{definition}{Definition}[section]

\def\bb{\Bbb}
\def\bpr{\begin{proof}}
\def\epr{\end{proof}}
\def\be{\begin{equation}}
\def\ee{\end{equation}}
\def\bea{\begin{eqnarray}}
\def\eea{\end{eqnarray}}
\def\bean{\begin{eqnarray*}}
\def\eean{\end{eqnarray*}}
\def\ba{\begin{abstract}}
\def\ea{\end{abstract}}
\def\bt{\begin{theorem}}
\def\et{\end{theorem}}
\def\bl{\begin{lemma}}
\def\el{\end{lemma}}
\def\bp{\begin{proposition}}
\def\ep{\end{proposition}}
\def\br{\begin{remark}}
\def\er{\end{remark}}
\def\bc{\begin{corollary}}
\def\ec{\end{corollary}}
\def\bd{\begin{definition}}
\def\ed{\end{definition}}
\def\non{\nonumber}

\def\bs{\bigskip}

\def\a{\alpha}
\def\b{\beta}

\begin{document}
 \baselineskip 20pt
\title{ The density theorem of a class of dilation-and-modulation systems on the half real line\thanks{Supported by the National Natural Science
Foundation of China (Grant No. 11271037).
}}
\author{Yun-Zhang Li$^\dagger$\,\,\,\,\ \ Ya-Hui Wang$^\ddagger$\\
College of Applied Sciences,\\ Beijing University of Technology,
Beijing 100124,  P. R. China\\E-mail:  $^\dagger$ yzlee@bjut.edu.cn, {\mbox{ }}$^\ddagger$ wangyahui@emails.bjut.edu.cn}\date{}
\maketitle{}
\begin{abstract} In the practice, time variable cannot be negative. The space $L^2(\bb R_+)$ of square integrable functions defined on the right half real line $\bb R_+$ models causal signal space. This paper focuses on a class of dilation-and-modulation systems in $L^2(\bb R_+)$. The density theorem for Gabor systems in $L^2(\bb R)$ states a necessary and sufficient condition for the existence of complete Gabor systems or Gabor frames in $L^2(\bb R)$ in terms of the index set alone-independently of window functions. The space $L^2(\bb R_+)$ admits no nontrivial Gabor system since $\bb R_+$ is not a group according to the usual addition. In this paper, we introduce a class of dilation-and-modulation systems in $L^2(\bb R_+)$ and the notion of $\Theta$-transform matrix. Using $\Theta$-transform matrix method we obtain the density theorem of the dilation-and-modulation systems in $L^2(\bb R_+)$ under the condition that $\log_ba$ is a positive rational number, where $a$ and $b$ are the dilation and modulation parameters respectively. Precisely, we prove that a necessary and sufficient condition for the existence of such a complete dilation-and-modulation system or dilation-and-modulation system frame in $L^2(\bb R_+)$ is that $\log_ba \leq 1$. Simultaneously, we obtain a $\Theta$-transform matrix-based expression of all complete dilation-and-modulation systems and all dilation-and-modulation system frames in $L^2(\bb R_+)$.
\end{abstract}

\noindent {\bf 2010 MS Classification}:  42C40, 42C15
\section{Introduction}
\setcounter{equation}{0}
A central part of
harmonic analysis deals with functions on groups and ways to decompose
such functions in terms of either series representations or integral representations
of certain ``basic functions".  Wavelet and Gabor frames are such  basic functions representing square integrable functions on $\bb R$.
 Given  a finite subset $\Psi$ of $L^{2}(\bb R)$,
Gabor frames of the form
 \be\label{gab}\{M_{mb}T_{na}\psi:\,m,\,n\in \bb Z,\,\psi\in \Psi\}\ee
 and  wavelet frames of the form
\be\label{wav}\{D_{a^{j}}T_{bk}\psi:\,j,\,k\in \bb Z, \,\psi\in \Psi\}\ee
  with $a$, $b>0$  have been extensively studied,
where the  translation
operator  $T_{x_{0}}$, the modulation operator $M_{x_{0}}$ with $x_{0}\in \bb R$ and the dilation operator $D_{c}$ with $0<c\ne 1$ are  defined by
$$T_{x_{0}}f(\cdot)=f(\cdot-x_{0}),{\mbox{ }}
M_{x_{0}}f(\cdot)=e^{2\pi ix_{0}\cdot}f(\cdot){\mbox{ and }}
D_{c}f(\cdot)=\sqrt{c}f(c\cdot)$$ for $f\in L^{2}(\bb R)$
respectively.  And Gabor analysis has been also extended to $l^2(\bb Z)$ and $\bb C^L$ ([\ref{O}, \ref{FS1}, \ref{FS2}, \ref{Gro01}, \ref{hei11}-\ref{her96}]).  Since $\bb R$, $\bb Z$ and $\bb C^L$  are all
locally compact abelian groups according to the usual addition and  topology, Gabor analysis in the
three cases has some similarity although there exist many differences in some aspects. The idea of considering frame theory on locally compact abelian groups has appeared in several publications including [\ref{bow15}, \ref{cab10}, \ref{O15}, \ref{dah94}, \ref{fei89}, \ref{kut06}, \ref{jakm16}, \ref{jakms16}].
Write $\bb R_+=(0,\,\infty)$. Given $a$, $b>0$, a measurable function $h$ defined on $\bb R_+$ is said to be {\it b-dilation periodic} if $h(b\cdot)=h(\cdot)$ on $\bb R_+$. Define  the sequence $\{\Lambda_{m}\}_{m\in\bb Z}$ of
$b$-dilation periodic functions by
\be \label{1111}
\Lambda_{m}(\cdot)=\frac{1}{\sqrt{b-1}}e^{\frac{2\pi im\cdot}{b-1}}\mbox{ on }[1,\,b)\mbox{ for each }m\in\bb Z.
\ee
Our focus in this paper will
be on  the {\it dilation-and-modulation systems} (${\cal MD}$-systems) in $L^2(\bb R_+)$ of the form:
\be \label{000}
{\cal MD}(\psi,\,a,\,b)=\left\{\,{\Lambda_{m}}D_{a^{j}}\psi:\,m,\,j\in \bb Z\,\right\}
\ee
with $\psi\in L^2(\bb R_+)$ under the following General setup:

\bs
{\bf General setup:}

(i) $a$ and $b$ are two  constants greater than 1.

(ii) $\log_{b}a=\frac{p}{q}$,  $p$ and $q$ are two coprime positive integers.

\bs
The space $L^2(\bb R_+)$ can be considered as a closed subspace of $L^2(\bb R)$ consisting of all
functions in $L^2(\bb R)$ which vanish outside $\bb R_+$. And it   models causal signal space. In the practice, time variable cannot be negative. Mathematically, we are inspired by the following observations to study  ${\cal MD}$-systems of the form (\ref{000}):

$\bullet$ The space $L^2(\bb R_+)$ admits
no nontrivial shift invariant system of the form $\{\,T_{na}\psi(\cdot):\,n\in\bb Z\,\}$ since $\psi=0$ is a unique $L^2(\bb R_+)$-solution to
$$T_{na}\psi(\cdot)=0{\mbox{ on }}(-\infty,\,0){\mbox{ for }}n\in\bb Z.$$
This implies that it  admits no nontrivial wavelet or Gabor system.

$\bullet$  An $\cal MD$-system in $L^2(\bb R_+)$ cannot be derived from  a wavelet system in the Hardy space
$H^{2}(\bb R)$ via  Fourier transform  since  the Fourier transform version of
$\{D_{a^{j}}T_{m}\psi:\,m,\,j\in \bb Z\}$ is
$\{M_{a^jm}D_{a^{j}}\hat\psi:\,m,\,j\in \bb Z\}
$.

$\bullet$ The space $L^2(\bb R_+)$ is not closed under the Fourier transform
since the Fourier transform of a compactly supported nonzero function in  $L^2(\bb R_+)$ lies outside this space.

$\bullet$ $\bb R$, $\bb Z$ and $\bb C^L$ are locally compact abelian groups according to the usual addition and topology, while $\bb R_+$ is not since the difference between two numbers in $\bb R_+$ may be negative. So the analysis in $L^2(\bb R_+)$ differs from that in $L^2(\bb R)$. Also $\bb R_+$ is a locally compact abelian group according to the usual multiplication and topology.

These observations inspire us to study multiplication-based frames for $L^2(\bb R_+)$ of the form (\ref{000}). In (\ref{000}), the dilation periodicity and expression on $[1,\,b)$ of $\Lambda_m$ are to match the dilation operation on $\psi$ and to apply the Fourier series theory.

The density theorem essentially states that
necessary and sufficient conditions for  a  Gabor system
to be complete, a frame, a Riesz basis, or a Riesz sequence in $L^2(\bb R)$ or general $L^2(\bb R^d)$ can be formulated in
terms of the index set alone-independently of the window function.
It has a long and very involved history from the one-dimensional rectangular lattice setting, to arbitrary lattices in higher dimensions, to irregular Gabor systems, and most recently
beyond the setting of Gabor systems to abstract localized systems. For details, see  [\ref{bag90}, \ref{bal06},  \ref{GHL13}, \ref{GL09}, \ref{hei07}, \ref{jak16}, \ref{jan13}, \ref{ly13}, \ref{rie81}] and references therein.

Let $a$ and $b$ be as in the general setup. We always write
\bea \label{2007} \beta=b^p~(\mbox{or }a^q).\eea
In this paper, we study the density theorem for the systems of the form (\ref{000}) in $L^2(\bb R_+)$. Section 2 is an auxiliary one. In this section, we introduce the notions of $\Theta_\beta$-transform and $\Theta_\beta$-transform matrix, and study their properties. In Section 3, using $\Theta_\beta$-transform matrix method we characterize complete ${\cal MD}$-systems and ${\cal MD}$-frames, obtain a parametrized expression of all complete systems and frames of the form (\ref{000}) in $L^2(\bb R_+)$, and finally we derive the density theorem. It turns out that $``\log_ba\leq 1"$ is necessary and sufficient for the existence of complete ${\cal MD}$-systems and ${\cal MD}$-frames (see Theorem \ref{Density}). Our main results are Theorems \ref{complete1}-\ref{Density}. They reduce constructing complete ${\cal MD}$-systems and ${\cal MD}$-frames to designing
$q\times p$ matrix-valued functions, and show that ${\cal MD}$-frames are significantly different from Gabor and wavelet frames for $L^2(\bb R)$. The lower bound $A$ and upper bound $B$ of an arbitrary ${\cal MD}$-frame are related by the inequality $a^{\frac{q-1}{p}}A\le B$, while it is not the case for a Gabor or wavelet frame. And an ${\cal MD}$-frame being a tight frame is equivalent to $``a=b"$ (see Remarks \ref{rem5} and \ref{rem6}).

\bs
Before proceeding, we introduce some notations and notions. Throughout this paper, the relation of quality, inclusion or inequality between two measurable sets is understood up to a set of measure zero, and the relation of quality or inequality between two measurable functions is understood in almost everywhere sense. We  denote by $I_t$ the $t\times t$ identity matrix, and by $\bb N_t$ the set
$$\bb N_t=\{0,\,1,\,2, \cdots,t-1\}$$
for $t\in\bb N$.

\bd\label{congr}
For a measurable set
$S\subset \bb R$, a collection  $\{S_i\}_{i\in \cal I}$ of measurable subsets of $S$
is called a partition of $S$ if $\chi_{_S}=\sum\limits_{i\in \cal I}\chi_{_{S_i}}$,
where $\chi_{_E}$ denotes the characteristic function of $E$ for a set $E$.
And for $\a>0$,  measurable  subsets $S_1$, $S_2$ of $\bb R$,  $T_1$, $T_2$  of $\bb R_+$ and a collection $\{T_i:\,i\in\cal I\}$ of  measurable subsets of $\bb R_+$ with $\cal I$ being at most countable,  we say that $S_1$ is {\it $\a\bb Z$-congruent} to $S_2$ if there exists a partition $\{S_{1,k}: \,k\in\bb Z\}$ of $S_1$ such that $\{S_{1,k}+\a k: k\in\bb Z\}$ is a partition of $S_2$,  that $T_1$ is {\it $\a^\bb Z$-dilation congruent} to $T_2$ if there exists a partition $\{T_{1,k}: k\in\bb Z\}$ of $T_1$ such that $\{\a^kT_{1,k}: \,k\in\bb Z\}$ is a partition of $T_2$,
 and that $\{T_i:\,i\in\cal I\}$ is {\it $\a^\bb Z$-dilation disjoint} if

$$T_i\cap\a^kT_j=\emptyset
$$
for $i$, $j\in\cal I$ and $k\in \bb Z$ with $(i,\,k)\neq(j,\,0)$.
\ed

\br\label{rem2} By Definition \ref{congr}, $S_2$ is {\it $\a\bb Z$-congruent} to $S_1$ $\rm{(}$$T_2$ is {\it $\a^\bb Z$-dilation congruent} to $T_1$$\rm{)}$ if $S_1$ is {\it $\a\bb Z$-congruent} to $S_2$ $\rm{(}$$T_1$ is {\it $\a^\bb Z$-dilation congruent} to $T_2$$\rm{)}$. So we say that $S_1$ and $S_2$ are {\it $\a\bb Z$-congruent} $\rm{(}$$T_1$ and $T_2$ are {\it $\a^\bb Z$-dilation congruent}$\rm{)}$ in this case. Also observe that $\bb Z$ is the superscript of $\alpha$ in the dilation congruence, and that only finitely many $S_{1,k}$ among $\{S_{1,k}: \,k\in\bb Z\}$ are nonempty if both $S_1$ and $S_2$ are bounded in addition. Similarly, only finitely many $T_{1,k}$ among $\{T_{1,k}: \,k\in\bb Z\}$ are nonempty if $T_1$ and $T_2$ are contained in some bounded subinterval $[M,\,N]$ of $\bb R_+$.\er

\br\label{rem3} By a standard argument, we see that, for $T_1,\,T_2\subset \bb R_+$, $T_1$ is {\it $\a^\bb Z$-dilation congruent} to $T_2$ if and only if $\log_\alpha T_1$ is {\it $\bb Z$-congruent} to $\log_\alpha T_2$, where
$$\log_\alpha T=\{\log_\alpha x:\,x\in T\} \mbox{ for } T\subset \bb R_+.$$
Similarly, $\{T_i:\,i\in\cal I\}$ is {\it $\a^\bb Z$-dilation disjoint} if and only if

$$\log_\alpha T_i\cap(k+\log_\alpha T_j)=\emptyset
$$
for $i$, $j\in\cal I$ and $k\in \bb Z$ with $(i,\,k)\neq(j,\,0)$.\er

Given $d,\,L\in \bb N$ and a measurable set $E\subset \bb R^d$, we denote by $L^2(E,\,\bb C^L)$ the vector-valued Hilbert space
$$L^2(E,\,\bb C^L)=\{f:\,\int_E\sum\limits_{l=1}^L|f_l(x)|^2dx<\infty\}$$
with the inner product
$$\langle f,\,g\rangle_{L^2(E,\,\bb C^L)}=\int_E\sum\limits_{l=1}^L f_l(x)\overline{g_l(x)}dx$$
for $f,\,g\in L^2(E,\,\bb C^L)$, where $f_l$ and $g_l$ denote the $l$-th components of $f$ and $g$ respectively. We denote by $L_{loc}^2(E,\,\bb C^L)$ the space of locally square integrable  vector-valued functions on $E$, i.e., the set of $f$ defined on $E$ satisfying
$$\int_F\sum\limits_{l=1}^L|f_l(x)|^2dx<\infty$$
for every bounded measurable subset $F$ of  $E$.

\bs
\section{ $\Theta_{\b}$-transform and $\Theta_{\b}$-transform matrix}
\setcounter{equation}{0}

Let $a$ and $b$ be as in the general setup, and $\beta$ be defined as in (\ref{2007}). This section is devoted to $\Theta_{\beta}$-transform matrix and related properties, which is an auxiliary one to following sections.

Define $\Theta_{\b}:\,L^2(\bb R_+)\to L^2_{loc}(\bb R_+\times\bb R)$ by
\be \label{2001}
\Theta_{\b}f(x,\xi)=\sum\limits_{l\in\bb Z}\b^{\frac{l}{2}}f(\b^{l}x)e^{-2\pi il\xi}
\ee
for $f\in L^{2}(\bb R_{+})$ and a.e. $(x,\,\xi)\in {\bb R_+}\times \bb R$.
It is well-defined due to
$$\int_{\b^j}^{\b^{j+1}}\sum\limits_{l\in\bb Z}\b^{l}|f(\b^{l}x)|^{2}dx=\|f\|_{L^{2}(\bb R_{+})}^{2}<\infty{\mbox{ for }}j\in\bb Z.$$
Using $\Theta_{\b}$ we define $\Gamma:\,L^2(\bb R_+)\to L^2_{loc}(\bb R_+\times\bb R,\,\bb C^p)$  by
\bea\Gamma f(x,\,\xi)=\left(b^{\frac{s}{2}}\Theta_\beta f(b^sx,\,\xi)\right)_{s\,\in \bb N_p}
\eea
for $f\in L^{2}(\bb R_{+})$ and a.e. $(x,\,\xi)\in {\bb R_+}\times \bb R$. Define the function sequence $\left\{e_{m,j}\right\}_{m,j\in\bb Z}$ by
\be \label{emj} e_{m,j}(x,\,\xi)=\Lambda_{m}(x)e^{2\pi ij\xi}{\mbox{ for }}(x,\,\xi)\in\bb R^2,\ee
where $\Lambda_{m}$ is  as in (\ref{1111}).

By a standard argument, we have the following two lemmas which  partially appeared in [\ref{lizh17}, Lemma 2.1].
\bl\label{onb}
{\rm (i) } $\{\Lambda_{m}:\,m\in\bb Z\}$ and $\{e_{m,j}:\,m,\,j\in\bb Z\}$ are orthonormal bases for $L^{2}([1,\,b))$ and $L^{2}([1,\,b)\times [0,\,1))$ respectively;

{\rm (ii) }
\be\label{2222}
\int_{[1,\,b)\times [0,\,1)}|f(x,\,\xi)|^{2}dxd\xi=\sum\limits_{j\in \bb Z}\sum\limits_{m\in \bb Z}\left|\int_{[1,\,b)\times [0,\,1)}f(x,\,\xi)\overline{e_{m,j}(x,\,\xi)}dxd\xi\right|^{2}
\ee for $f\in L^{1}([1,\,b)\times [0,\,1))$.
\el

\bl \label{map}

{\rm (i) } $\Theta_{\beta}$ has the quasi-periodicity:
$$\Theta_{\beta}f(\beta^{j}x,\,\xi+m)=\beta^{-\frac{j}{2}}e^{2\pi ij\xi}\Theta_{\beta}f(x,\,\xi) \mbox{ for }f \in L^2(\bb R_+),\,m,\,j\in \bb Z \mbox{ and a.e. }(x,\,\xi)\in \bb R_+\times \bb R.$$

{\rm (ii) }For $f \in L^2(\bb R_+)$ and $(j,\,m,\,r) \in \bb Z\times\bb Z\times\bb N_q$,
$$\Theta_{\beta}\Lambda_{m}D_{a^{jq+r}}f(x,\,\xi)=a^{\frac{r}{2}}e_{m,\,j}(x,\,\xi)\Theta_{\beta}f(a^rx,\,\xi),$$
$$\Gamma\Lambda_{m}D_{a^{jq+r}}f(x,\,\xi)=e_{m,\,j}(x,\,\xi)\left(a^{\frac{r}{2}}b^{\frac{s}{2}}
\Theta_{\beta}f(a^rb^sx,\,\xi)\right)_{s\in \bb N_p}$$
for a.e. $(x,\,\xi)\in \bb R_+\times \bb R$.

{\rm (iii) } The mappings $\Theta_{\beta}$ and $\Gamma$ are unitary operators from  $L^{2}(\bb R_{+})$ onto $L^{2}([1,\,\beta)\times[0,\,1))$ and $L^{2}([1,\,b)\times[0,\,1),\,\bb C^p)$ respectively.
\el

\br\label{rem1} By Lemma \ref{map} (i) and (iii), a function $f\in L^2(\bb R_+)$ is uniquely determined by $\Theta_\beta f$ restricted on $S\times [0,\,1)$ with $S$ being a set which is {\it $\beta^\bb Z$}-dilation congruent to $[1,\,\beta)$. Let $S$ be a set {\it $\beta^\bb Z$}-dilation congruent to $[1,\,\beta)$, and $\widetilde{f}\in L^2\left(S\times [0,\,1)\right)$. Define $f$ on $\bb R_+$ by
$$f(\beta^jx)=\beta^{-\frac{j}{2}}\int_0^1\widetilde{f}(x,\,\xi)e^{2\pi ij\xi}d\xi \mbox{ for }j\in \bb Z \mbox{ and a.e. }x\in S.$$
Then $f$ is well-defined and the unique function satisfying $\Theta_\beta f=\widetilde{f}$ a.e. on $S\times [0,\,1)$. Indeed, since
$$\|f\|_{L^2(\bb R_+)}^2=\int_S\sum\limits_{j\in \bb Z}\beta^j|f(\beta^jx)|^2dx,$$
we have
$$\|f\|_{L^2(\bb R_+)}^2=\int_{S\times [0,\,1)}\left|\widetilde{f}(x,\,\xi)\right|^2dxd\xi<\infty.$$
Also observing that $\widetilde{f}(x,\,\cdot)\in L^2[0,\,1)$ for a.e. $x\in S$ leads to
$$\widetilde{f}(x,\,\xi)=\sum\limits_{j\in \bb Z}\left(\int_0^1\widetilde{f}(x,\,\xi)e^{2\pi ij\xi}d\xi\right)e^{-2\pi ij\xi}$$
for a.e. $(x,\,\xi)\in S\times [0,\,1)$. This is in turn equivalent to
$$\Theta_\beta f=\widetilde{f} \mbox{ a.e. on }S\times [0,\,1)$$
by the definition of $f$.\er

\bl \label{congruent2} Let $a$ and $b$ be as in the general setup. Then the collection $\{a^rb^s[1,\,a^{\frac{1}{p}}):\,(r,\,s)\in \bb N_q\times \bb N_p\}$ is {\it $\beta^\bb Z$}-dilation disjoint, and the set $\bigcup\limits_{r\in \bb N_q}\bigcup\limits_{s\in \bb N_p}a^rb^s[1,\,a^{\frac{1}{p}})$ is {\it $\beta^\bb Z$}-dilation congruent to $[1,\,\beta)$.
\el
\bpr By Remark \ref{rem3}, the conclusion of this lemma is equivalent to the fact that
\bea\label{equi}\left(\frac{r_1}{q}+\frac{s_1}{p}+[0,\,\frac{1}{pq})\right)\bigcap\left(k+\frac{r_2}{q}+\frac{s_2}{p}+[0,\,\frac{1}{pq})\right)
=\emptyset\eea
if $(r_1,\,s_1,\,k)\neq(r_2,\,s_2,\,0)$ for $(r_1,\,s_1)$, $(r_2,\,s_2)\in \bb N_q\times \bb N_p$ and $k\in \bb Z$, and that the set $$\bigcup\limits_{r\in \bb N_q}\bigcup\limits_{s\in \bb N_p}\left(\frac{r}{q}+\frac{s}{p}+[0,\,\frac{1}{pq})\right)$$ is {\it $\bb Z$}-congruent to $[0,\,1)$.

Obviously, (\ref{equi}) can be rewritten as
\bea\label{equi1}\left(\frac{pr_1+qs_1}{pq}+[0,\,\frac{1}{pq})\right)\bigcap\left(\frac{kpq+pr_2+qs_2}{pq}+[0,\,\frac{1}{pq})\right)
=\emptyset.\eea
It is easy to check that, if $(r_1,\,s_1,\,k)\neq(r_2,\,s_2,\,0)$, then
$$pr_1+qs_1\neq kpq+pr_2+qs_2,$$
and thus (\ref{equi1}) holds. This leads to (\ref{equi}). Also observe that
$$\sum\limits_{r\in \bb N_q}\sum\limits_{s\in \bb N_p}\left|\frac{r}{q}+\frac{s}{p}+[0,\,\frac{1}{pq})\right|=1.$$
It follows that
\bea\label{equi2}\bigcup\limits_{r\in \bb N_q}\bigcup\limits_{s\in \bb N_p}\left(\frac{r}{q}+\frac{s}{p}+[0,\,\frac{1}{pq})\right)\eea
is {\it $\bb Z$}-congruent to $[0,\,1)$.  The proof is completed.
\epr

\bd\label{trma} Let $a$ and $b$ be as in the general setup. We associate every $\psi\in L^2(\bb R_+)$ with a $q\times p$ matrix-valued function $\Psi$ $\rm{(}$also called $\Theta_\beta$-transform matrix$\rm{)}$ on $\bb R_+\times\bb R$ by
$$\Psi(x,\,\xi)=\left(\Psi_{r,\,s}(x,\,\xi)\right)_{r\in \bb N_q,\,s\in \bb N_p}$$
with
$$\Psi_{r,\,s}(x,\,\xi)=a^{\frac{r}{2}}b^{\frac{s}{2}}\Theta_\beta\psi(a^rb^sx,\,\xi)$$
for a.e. $(x,\,\xi)\in \bb R_+\times \bb R$.
\ed

\br\label{rem4}By Remark \ref{rem1} and Lemma \ref{congruent2}, a function $\psi\in L^2(\bb R_+)$ is uniquely determined by the values of its $\Theta_\beta$-transform matrix $\Psi$ on $[1,\,a^{\frac{1}{p}})\times[0,\,1)$. For a $q\times p$ matrix-valued function $\Psi$ on $[1,\,a^{\frac{1}{p}})\times[0,\,1)$ with $L^2([1,\,a^{\frac{1}{p}})\times[0,\,1))$-entries, there exists a unique $\psi\in L^2(\bb R_+)$ such that
$$\left(a^{\frac{r}{2}}b^{\frac{s}{2}}\Theta_\beta\psi(a^rb^sx,\,\xi)\right)_{r\in \bb N_q,\,s\in \bb N_p}=\Psi(x,\,\xi) \mbox{ for }(x,\,\xi)\in [1,\,a^{\frac{1}{p}})\times[0,\,1).$$
And the unique function $\psi$ is defined by
$$\psi(\beta^ja^rb^sx)=\beta^{-\frac{j}{2}}a^{-\frac{r}{2}}b^{-\frac{s}{2}}\int_0^1\Psi_{r,s}(x,\,\xi)e^{2\pi ij\xi}d\xi \mbox{ for }j\in \bb Z,\,(r,\,s)\in \bb N_q\times \bb N_p \mbox{ and a.e. }x\in [1,\,a^{\frac{1}{p}}).$$\er

\bl \label{frame3}  Let $p$ and $q$ be two coprime positive integers, and $p,\,q>1$. Then there exists a unique $(r,\,s)\in (\bb N_q\backslash\{0\}) \times \left(\bb N_p\backslash\{0\}\right)$ such that
\bea\label{equa2}pr+qs=pq+1.\eea\el
\bpr First, we prove the existence of $(r,\,s)$ satisfying (\ref{equa2}). Without loss of generality, we assume that $p<q$. Since $p$ and $q$ are coprime, there exists $(s_1,\,\mu_1)\in \bb Z^2$ such that
\bea\label{equa3}1=q s_1+p\mu_1.\eea
Also observe that $s_1$ has the decomposition $$s_1=mp+s \mbox{ for some }(m,\,s)\in \bb Z\times \bb N_p.$$
It follows that
\bea\label{equa4}1=q s-\mu p \mbox{ with } \mu=-m q-\mu_1.\eea
Take $r=q-\mu$. Then
\bea\label{equa5}pr+qs=pq-\mu p+qs=pq+1\eea
by (\ref{equa4}). If $rs=0$, then $pr+qs<pq$, contradicting (\ref{equa5}). From (\ref{equa4}), we have
$$r=q-\mu=q-\frac{qs-1}{p}.$$
It follows that $r<q$ due to $s\neq0$, and
$$r\geq q-\frac{q(p-1)-1}{p}=\frac{q+1}{p}>1$$
due to $s\leq p-1$ and $p<q$. Therefore,
$(r,\,s)\in (\bb N _q\backslash\{0\})\times(\bb N_p\backslash\{0\})$, and $(r,\,s)$ satisfies (\ref{equa2}).

Next we prove the uniqueness of $(r,\,s)$ in (\ref{equa2}). Suppose that $(r_1,\,s_1)$, $(r_2,\,s_2)\in (\bb N _q\backslash\{0\})\times(\bb N_p\backslash\{0\})$ satisfy
$$pr_1+qs_1=pq+1,~pr_2+qs_2=pq+1.$$
Then $pr_1+qs_1=pr_2+qs_2$, and thus
$$p(r_1-r_2)=q(s_2-s_1).$$
This implies that $p\big|(s_2-s_1)$ and $q\big|(r_1-r_2)$ since $p$ and $q$ are coprime. It follows that
$$(r_1,\,s_1)=(r_2,\,s_2)$$
due to $r_1,\,r_2\in \bb N_q$ and $s_1,\,s_2\in \bb N_p$. The proof is completed.
\epr

Under the hypothesis of Lemma \ref{frame3}, suppose that $(r',\,s')\in \bb N_q\times \bb N_p$ is such that $pr'+qs'=pq+1$. Define
\bea\label{matr}{\cal L}_q(\xi)=\begin{pmatrix} 0 & I_{q-r'} \\ e^{2\pi i\xi}I_{r'} & 0 \end{pmatrix} \mbox{ and }{\cal R}_p(\xi)=\begin{pmatrix} 0 & I_{s'} \\ e^{-2\pi i\xi}I_{p-s'} & 0 \end{pmatrix}.\eea
Then they are uniquely determined by $p$ and $q$.

\bl  \label{frame4} Let $a$ and $b$ be as in the general setup, and $\psi\in L^2(\bb R_+)$. We have

{\rm (i)} $$\Psi(a^{lq+m}x,\,\xi)=a^{-\frac{lq+m}{2}}(e^{2\pi i\xi})^lU_m(\xi)\Psi(x,\xi)$$ for $(l,\,m)\in \bb Z\times\bb N_q$ and a.e. $(x, \,\xi)\in \bb R_+\times \bb R$,
where $U_m(\xi)=\begin{pmatrix} 0 & I_{q-m} \\ e^{2\pi i\xi}I_m & 0 \end{pmatrix}.$

{\rm (ii)} If $p,\,q>1$, then, for a.e. $(x, \,\xi)\in \bb R_+\times \bb R$,
\bea \Psi(a^{\frac{1}{p}}x,\xi)=a^{-\frac{1}{2p}}{\cal L}_q(\xi)\Psi(x,\xi){\cal R}_p(\xi).\nonumber\eea
\el
\bpr (i) For $l\in \bb N$, $(r,\,s)\in \bb N_q\times \bb N_p$ and $(x,\,\xi)\in \bb R_+\times \bb R$, we have
\begin{align}\label{equi3}
\Psi_{r,\,s}(a^{lq}x,\,\xi)&=a^{\frac{r}{2}}b^{\frac{s}{2}}\Theta_\beta\psi(a^rb^s\beta^lx,\,\xi)\non\\
&=a^{-\frac{lq}{2}}e^{2\pi i l\xi}\Psi_{r,\,s}(x,\,\xi)
\end{align}
by Lemma \ref{map} (i). This leads to
$$\Psi_{r,\,s}(x,\,\xi)=a^{-\frac{lq}{2}}e^{2\pi i l\xi}\Psi_{r,\,s}(a^{-lq}x,\,\xi)$$
by substituting $a^{lq}x$ for $x$ in (\ref{equi3}), equivalently,
$$\Psi_{r,\,s}(a^{-lq}x,\,\xi)=a^{\frac{lq}{2}}e^{-2\pi i l\xi}\Psi_{r,\,s}(x,\,\xi).$$
Together with (\ref{equi3}), it follows that
\bea\label{equi4}\Psi_{r,\,s}(a^{lq}x,\,\xi)=a^{-\frac{lq}{2}}e^{2\pi i l\xi}\Psi_{r,\,s}(x,\,\xi)\eea
for $l\in \bb Z$, $(r,\,s)\in \bb N_q\times \bb N_p$ and a.e. $(x,\,\xi)\in \bb R_+\times \bb R.$

For $m\in \bb N_q$, $(r,\,s)\in \bb N_q\times \bb N_p$ and $(x,\,\xi)\in \bb R_+\times \bb R$, we have
\bea\label{equi5}\Psi_{r,\,s}(a^mx,\,\xi)=\left\{\begin{array}{ll}
a^{-\frac{m}{2}}\Psi_{r+m,\,s}(x,\,\xi) & \hbox{${\text{ if }}0\leq r\leq q-1-m$};\\
a^{-\frac{m}{2}}e^{2\pi i\xi}\Psi_{r+m-q,\,s}(x,\,\xi) & \hbox{${\text{ if }}q-m\leq r\leq q-1$}
\end{array} \right.\eea
by Lemma \ref{map} (i). So
\begin{align*}
\Psi_{r,\,s}(a^{lq+m}x,\,\xi)&=a^{-\frac{lq}{2}}e^{2\pi i l\xi}\Psi_{r,\,s}(a^mx,\,\xi)\\
&=\left\{
\begin{aligned}
&a^{-\frac{lq+m}{2}}\Psi_{r+m,\,s}(x,\,\xi)~~~~~~~~~~~\text{ if } 0\leq r\leq q-1-m; \\
&a^{-\frac{lq+m}{2}}e^{2\pi il\xi}\Psi_{r+m-q,\,s}(x,\,\xi)~\text{ if } q-m\leq r\leq q-1
\end{aligned}
\right.\end{align*}
for $(l,\,m)\in \bb Z\times \bb N_q$, $(r,\,s)\in \bb N_q\times \bb N_p$ and a.e. $(x,\,\xi)\in \bb R_+\times \bb R$ by (\ref{equi4}) and (\ref{equi5}). This leads to (i).

(ii) Write $$\Psi(x,\,\xi)=\begin{pmatrix}{\cal B}_{0,0}(x,\,\xi) & {\cal B}_{0,1}(x,\,\xi)\\{\cal B}_{1,0}(x,\,\xi) & {\cal B}_{1,1}(x,\,\xi)\end{pmatrix}$$
for a.e. $(x,\,\xi)\in \bb R_+\times \bb R$, where ${\cal B}_{0,0}$ is of order $r'\times s'$.

Since $pr'+qs'=pq+1$, we have
$$\beta^{-1}a^{r'}b^{s'}=a^{\frac{1}{p}}.$$
It follows that
$$\Psi_{r,\,s}(a^{\frac{1}{p}}x,\,\xi)=a^{\frac{r}{2}}b^{\frac{s}{2}}\Theta_\beta\psi(a^{r+r'}b^{s+s'}\beta^{-1}x,\,\xi),$$
and thus
\begin{align*}
\Psi_{r,\,s}(a^{\frac{1}{p}}x,\,\xi)&=\beta^{\frac{1}{2}}e^{-2\pi i\xi}a^{\frac{r}{2}}b^{\frac{s}{2}}\Theta_\beta\psi(a^{r+r'}b^{s+s'}x,\,\xi)\\
&=\left\{
\begin{aligned}
&a^{-\frac{1}{2p}}e^{-2\pi i\xi}\Psi_{r+r',\,s+s'}(x,\,\xi)~~~~~\text{ if } r'\leq r+r'\leq q-1,\,s'\leq s+s'\leq p-1; \\
&a^{-\frac{1}{2p}}\Psi_{r+r',\,s+s'-p}(x,\,\xi)~~~~~~~~~~\text{ if } r'\leq r+r'\leq q-1,\,s+s'\geq p;\\
&a^{-\frac{1}{2p}}\Psi_{r+r'-q,\,s+s'}(x,\,\xi)~~~~~~~~~~\text{ if } r+r'\geq q,\,s'\leq s+s'\leq p-1;\\
&a^{-\frac{1}{2p}}e^{2\pi i\xi}\Psi_{r+r'-q,\,s+s'-p}(x,\,\xi)~ \text{ if } r+r'\geq q,\,s+s'\geq p
\end{aligned}
\right.
\end{align*}
for $(r,\,s)\in \bb N_q\times \bb N_p$ and a.e. $(x,\,\xi)\in \bb R_+\times \bb R$ by Lemma \ref{map} (i) and a simple computation. So
$$\Psi(a^{\frac{1}{p}}x,\,\xi)=a^{-\frac{1}{2p}}\begin{pmatrix}e^{-2\pi i\xi}{\cal B}_{1,1}(x,\,\xi) & {\cal B}_{1,0}(x,\,\xi)\\{\cal B}_{0,1}(x,\,\xi) & e^{2\pi i\xi}{\cal B}_{0,0}(x,\,\xi)\end{pmatrix}$$
for a.e. $(x,\,\xi)\in \bb R_+\times \bb R$. This implies (ii). The proof is completed.
\epr

\bl \label{3333}Let $a$ and $b$ be as in the general setup. Then, for any $A,\,B>0$,
\bea \label{444}AI_p\leq \Psi^{\ast}(x,\,\xi)\Psi(x,\,\xi) \leq BI_p \mbox{ for a.e. }(x,\,\xi)\in\,[1,\,b)\times[0,\,1),\eea
if and only if
\bea \label{4444}a^{\frac{q-1}{p}}AI_p\leq \Psi^{\ast}(x,\,\xi)\Psi(x,\,\xi) \leq BI_p{\mbox{ for a.e. }}
(x,\,\xi)\in [1,\,a^{\frac{1}{p}})\times[0,\,1).\eea
\el
\bpr When $q=1$, we have $b=a^\frac{1}{p}$, and thus (\ref{444}) is exactly (\ref{4444}). Next we prove their equivalence for the case $q>1$.

By Lemma \ref{frame4} and a simple computation, we have
$$\Psi^{\ast}(a^{\frac{l}{p}}x,\,\xi)\Psi(a^{\frac{l}{p}}x,\,\xi)=\left\{
\begin{aligned}
&a^{-l}\Psi^{\ast}(x,\,\xi)\Psi(x,\,\xi)~~~~~~~~~~~~~~~~~~~~~~~~~~~\text{if } p=1; \\
&a^{-\frac{l}{p}}\left({\cal R}_p^{\ast}(\xi)\right)^l\Psi^{\ast}(x,\,\xi)\Psi(x,\,\xi)\left({\cal R}_p(\xi)\right)^l~~\text{if } p>1
\end{aligned}
\right.$$
for $l\in \bb N_q$ and a.e. $(x,\,\xi)\in [1,\,a^{\frac{1}{p}})\times[0,\,1)$, where ${\cal R}_p(\xi)$ is as in (\ref{matr}). Also observe that
$$[1,\,b)=\bigcup\limits_{l\in \bb N_q}a^{\frac{l}{p}}[1,\,a^{\frac{1}{p}})$$
and that ${\cal R}_p(\xi)$ is unitary. It follows that (\ref{444}) holds if and only if
 $$AI_p\leq a^{-\frac{l}{p}}\Psi^{\ast}(x,\,\xi)\Psi(x,\,\xi) \leq BI_p \mbox{ for } l\in \bb N_q \mbox{ and a.e. }
(x,\,\xi)\in [1,\,a^{\frac{1}{p}})\times[0,\,1),$$
which is in turn equivalent to (\ref{4444}). The proof is completed.
\epr

\bs
\section{The density theorem}
\setcounter{equation}{0}

Let $a$ and $b$ be as in the general setup. In this section, using $\Theta_\beta$-transform matrix method, we characterize complete ${\cal MD}$-systems and ${\cal MD}$-frames, present a parameterized expression of them, and derive the density theorem for ${\cal MD}$-systems of the form ${\cal MD}(\psi,\,a,\,b)$ in $L^2(\bb R_+)$.

\bl\label{unita0} Let $a$ and $b$ be as in the general setup, and $\psi \in L^2(\bb R_+)$. Then
\be \label{202}
\sum\limits_{m,\,j\in \bb Z}\left|\langle f,\,{\Lambda_{m}}D_{a^{j}}\psi\rangle_{L^{2}(\bb R_{+})}\right|^{2}=\int_{[1,b)\times[0,1)}\left\| \overline{\Psi(x,\,\xi)}\Gamma f(x,\xi)\right\|_{\bb C^q}^{2}dxd\xi
\ee
for $f\in L^2(\bb R_+)$.
\el
\bpr By Lemma \ref{map} (ii) and (iii), we have
\begin{align*}
\langle f,\,{\Lambda_{m}}D_{a^{jq+r}}\psi\rangle_{L^2(\bb R_+)}
&=\langle\Gamma f,\,\Gamma{\Lambda_{m}}D_{a^{jq+r}}\psi \rangle_{L^2([1,\,b)\times[0,\,1),\,\bb C^p)}\\
&=\int_{[1,\,b)\times[0,\,1)}\left(\overline{\Psi(x,\,\xi)}\Gamma f(x,\,\xi)\right)_r\overline{e_{m,\,j}(x,\,\xi)}dxd\xi
\end{align*}
for $f\in L^2(\bb R_+)$ and $(r,\,j)\in \bb N_q\times\bb Z$. It follows that
\begin{align*}
\sum\limits_{m,\,j\in \bb Z}\left|\langle f,\,{\Lambda_{m}}D_{a^{j}}\psi\rangle_{L^{2}(\bb R_{+})}\right|^{2}&
=\sum\limits_{r=0}^{q-1}\sum\limits_{j\in \bb Z}\sum\limits_{m\in \bb Z}\left|\langle  f,\,{\Lambda_{m}}D_{a^{jq+r}}\psi\rangle_{L^{2}(\bb R_{+})}\right|^{2}\\&
=\int_{[1,\,b)\times[0,\,1)}\left\|\overline{\Psi(x,\,\xi)}\Gamma f(x,\,\xi)\right\|_{\bb C^q}^2dxd\xi
\end{align*}
by Lemma \ref{onb} (ii). The proof is completed.
\epr

The following lemma is borrowed from \cite{hb17}, and it is a variation of [\ref{hb17}, Corollary 2.4].

\bl \label{matrixexp} An arbitrary $\mu\times \nu$ matrix-valued measurable function ${\cal A}$ on a measurable set $E$ in $\bb R^d$ must have the form
$${\cal A}(\cdot)=U(\cdot)\begin{pmatrix}\Lambda(\cdot) & 0\\ 0 & 0\end{pmatrix}V(\cdot)  \mbox{ a.e. on }E,$$
where $U(\cdot)$ and $V(\cdot)$ are $\mu\times \mu$ and $\nu\times \nu$ unitary matrix-valued measurable functions on $E$ respectively, and $\Lambda(\cdot)$ is a min$(\mu,\,\nu)\times$min$(\mu,\,\nu)$ diagonal matrix-valued measurable function on $E$.
\el

By an easy application of the spectrum theorem for self-adjoint matrices (see also [\ref{Dau}, p.978]), we have
 \bl\label{per} Let ${\cal A}(\cdot)$ be an $m\times n$ matrix-valued measurable function defined on a measurable set $E$. Then the orthogonal projection operator $P_{ker({\cal A}(\cdot))}$ from $\bb C^n$ onto the kernel space $ker({\cal A}(\cdot))$ of ${\cal A}(\cdot)$ is measurable on $E$, and
 $$P_{ker({\cal A}(\cdot))}=\lim\limits_{n\rightarrow \infty}exp(-n{\cal A}^{\ast}(\cdot){\cal A}(\cdot)).$$
 \el
\bt\label{complete1}Let $a$ and $b$ be as in the general setup. Then, for $\psi\in L^{2}(\bb R_{+})$, ${\cal MD}(\psi,\,a,\,b)$ is complete in $L^{2}(\bb R_{+})$ if and only if
\be \label{2005}
rank\left(\Psi(x,\xi)\right)=p\,\,\,{\mbox{ for a.e. }}(x,\,\xi)\in [1,\,a^{\frac{1}{p}})\times[0,\,1).
\ee\et
\bpr The system ${\cal MD}(\psi,\,a,\,b)$ is complete in $L^2(\bb R_+)$ if and only if $f=0$ is a unique solution to
\bea\label{equa9}\sum\limits_{m,\,j\in \bb Z}\left|\langle f,\,\Lambda_{m}D_{a^{j}}\psi\rangle_{L^{2}(\bb R_{+})}\right|^{2}=0\eea
in $L^2(\bb R_+)$, equivalently, $f=0$ is a unique solution to
\bea\label{equa10}\overline{\Psi(x,\,\xi)}\Gamma f(x,\xi)=0 \mbox{ for a.e. }(x,\,\xi)\in [1,\,b)\times[0,\,1)\eea
in $L^2(\bb R_+)$ by Lemma \ref{unita0}. By Lemma \ref{map} (iii), this is also equivalent to $F=0$ is a unique solution to
\bea\label{thc1}\overline{\Psi(x,\,\xi)}F(x,\xi)=0 \mbox{ for a.e. }(x,\,\xi)\in [1,\,b)\times[0,\,1)\eea
in $L^2([1,\,b)\times[0,\,1),\,\bb C^p)$. Also by Lemma \ref{frame4}, (\ref{2005}) holds if and only if
\bea\label{thc2} rank(\Psi(x,\,\xi))=p \mbox{ for a.e. }(x,\,\xi)\in [1,\,b)\times[0,\,1).\eea
So it is enough to prove that (\ref{thc2}) is equivalent to $F=0$ being a unique solution to (\ref{thc1}) in $L^2([1,\,b)\times[0,\,1),\,\bb C^p)$. Obviously, (\ref{thc2}) implies $F=0$ being a unique solution to (\ref{thc1}) in $L^2([1,\,b)\times[0,\,1),\,\bb C^p)$. Next we prove that (\ref{thc1}) has a nonzero solution in $L^2([1,\,b)\times[0,\,1),\,\bb C^p)$ if (\ref{thc2}) does not hold. Suppose (\ref{thc2}) does not hold. Then there exists some $E\subset[1,\,b)\times[0,\,1)$ with $|E|>0$ such that
\bea\label{thc3} rank(\Psi(x,\,\xi))<p \mbox{ for }(x,\,\xi)\in E.\eea
Let $P_{ker(\overline{\Psi(x,\,\xi)})}$ be the orthogonal projection from $\bb C^p$ onto the kernel space $ker(\overline{\Psi(x,\,\xi)})$ of $\overline{\Psi(x,\,\xi)}$, and $\{e_1,\,e_2,\,\cdots,\,e_p\}$ be the canonical orthonormal basis for $\bb C^p$, i.e., every $e_l$ with $1\le l\le p$ is the vector in $\bb C^p$ with the $l$-th component being 1 and others being zero. Then
$$span\{P_{ker(\overline{\Psi(x,\,\xi)})}e_l:\,1\le l \le p\}=ker(\overline{\Psi(x,\,\xi)}).$$
Observe that $ker(\overline{\Psi(x,\,\xi)})\neq \{0\}$ for $(x,\,\xi)\in E$ by (\ref{thc3}). It follows that there exist $1\le l_0\le p$ and $E'\subset E$ with $|E'|>0$ such that
$$P_{ker(\overline{\Psi(x,\,\xi)})}e_{l_0}\neq 0 \mbox{ for }(x,\,\xi)\in E'.$$
Take $F(x,\,\xi)=\chi_{_{E'}}(x,\,\xi)
P_{ker(\overline{\Psi(x,\,\xi)})}e_{l_0}$ for $(x,\,\xi)\in [1,\,b)\times[0,\,1)$. Then $F$ is measurable, nonzero and $\|F(x,\,\xi)\|_{\bb C^p}\le 1$ by Lemma \ref{per}. It is obvious that $F$ is a nonzero solution to (\ref{thc1}) in $L^2([1,\,b)\times[0,\,1),\,\bb C^p)$. The proof is completed.
\epr

\br\label{new} Observe that $\Psi$ is a $q\times p$ matrix. Theorem \ref{complete1} shows that the inequality
$$p\le q~\rm{( equivalently},~\log_ba\le 1\rm{)}$$
is necessary for the existence of complete ${\cal MD}$-systems in $L^2(\bb R_+)$.\er

\bt \label{frame1} Let $a$ and $b$ be as in the general setup. Then, for $\psi\in L^{2}(\bb R_{+})$, ${\cal MD}(\psi,\,a,\,b)$ is a frame for $L^2(\bb R_+)$ with frame bounds $A$ and $B$ if and only if
$$a^{\frac{q-1}{p}}AI_p\leq\Psi^{\ast}(x,\,\xi)\Psi(x,\xi)\leq BI_p \mbox{ for a.e. }(x,\,\xi)\in [1,\,a^{\frac{1}{p}})\times[0,\,1).$$\et
\bpr By Lemma \ref{unita0} and Lemma \ref{map} (iii), ${\cal MD}(\psi,\,a,\,b)$ is a frame for $L^2(\bb R_+)$ with frame bounds $A$ and $B$ if and only if
\begin{align}\label{equa12}A\int_{[1,\,b)\times[0,\,1)}\left\|\Gamma f(x,\,\xi)\right\|_{\bb C^p}^2dxd\xi&\leq\int_{[1,\,b)\times[0,\,1)}\left\|\overline{\Psi(x,\,\xi)}\Gamma f(x,\,\xi)\right\|_{\bb C^q}^2dxd\xi\non\\
&\leq B\int_{[1,\,b)\times[0,\,1)}\left\|\Gamma f(x,\,\xi)\right\|_{\bb C^p}^2dxd\xi.\end{align}
By a standard argument, (\ref{equa12}) is equivalent to
\bea\label{equa13}AI_p\leq\Psi^{\ast}(x,\,\xi)\Psi(x,\,\xi)\leq BI_p \mbox{ for a.e. }(x,\,\xi)\in [1,\,b)\times[0,\,1).\eea
This is equivalent to
\bea a^{\frac{q-1}{p}}AI_p\leq\Psi^{\ast}(x,\,\xi)\Psi(x,\,\xi)\leq BI_p \mbox{ for a.e. }(x,\,\xi)\in [1,\,a^{\frac{1}{p}})\times[0,\,1)\eea
by Lemma \ref{3333}.
\epr

By Remark \ref{new}, the inequality $\log_ba\le 1$ is necessary for the existence of a complete ${\cal MD}$-system ${\cal MD}(\psi,\,a,\,b)$ in $L^2(\bb R_+)$. The following theorem shows that this condition is also sufficient for the existence of a complete system ${\cal MD}(\psi,\,a,\,b)$ in $L^2(\bb R_+)$. Simultaneously, it shows that this condition is also sufficient for the existence of an ${\cal MD}$-frame ${\cal MD}(\psi,\,a,\,b)$ for $L^2(\bb R_+)$. Again by Remark \ref{rem4} and Lemma \ref{matrixexp}, it provides a parameterized expression of all complete ${\cal MD}$-systems and all ${\cal MD}$-frames in $L^2(\bb R_+)$ of the form ${\cal MD}(\psi,\,a,\,b)$.

\bt \label{Density1} Let $a$ and $b$ be as in the general setup, and $\log_ba\leq1$. Assume that $\lambda_s(\cdot,\,\cdot)\in L^2([1,\,a^{\frac{1}{p}})\times[0,\,1))$ with $s\in \bb N_p$. Define
$$\Lambda(\cdot,\,\cdot)=diag\left(\lambda_0(\cdot,\,\cdot),\,\lambda_1(\cdot,\,\cdot),\,\cdots,\,\lambda_{p-1}(\cdot,\,\cdot)\right)$$
and $\psi$ by a $q\times p$ matrix-valued function
$$\Psi(\cdot,\,\cdot)=U(\cdot,\,\cdot)\begin{pmatrix}\Lambda(\cdot,\,\cdot)\\0\end{pmatrix}V(\cdot,\,\cdot)$$
on $[1,\,a^{\frac{1}{p}})\times[0,\,1)$, where $U(\cdot,\,\cdot)$ and $V(\cdot,\,\cdot)$ are $q\times q$ and $p\times p$ unitary matrix-valued measurable functions defined on $[1,\,a^{\frac{1}{p}})\times[0,\,1)$ respectively. Then

{\rm (i) }${\cal MD}(\psi,\,a,\,b)$ is complete in $L^2(\bb R_+)$ if and only if
$$\lambda_0(\cdot,\,\cdot)\lambda_1(\cdot,\,\cdot)\cdots\lambda_{p-1}(\cdot,\,\cdot)\neq 0 \mbox{ a.e. on }[1,\,a^{\frac{1}{p}})\times[0,\,1).$$

{\rm (ii) }${\cal MD}(\psi,\,a,\,b)$ is a frame for $L^2(\bb R_+)$ if and only if there exist constants $0<A\leq B<\infty$ such that
$$A\leq \left|\lambda_s(\cdot,\,\cdot)\right|\leq B \mbox{ for }s\in \bb N_p \mbox{ a.e. on }[1,\,a^{\frac{1}{p}})\times[0,\,1).$$
\et
\bpr (i) Since rank$\left(\Psi\right)$=rank$\left(\Lambda\right)$, we have
$$\mbox{rank}\left(\Psi(x,\,\xi)\right)=p \mbox{ for a.e. }(x,\,\xi)\in [1,\,a^{\frac{1}{p}})\times[0,\,1)$$
if and only if $\lambda_0(x,\,\xi)\lambda_1(x,\,\xi)\cdots\lambda_{p-1}(x,\,\xi)\neq 0$ for a.e. $(x,\,\xi)\in [1,\,a^{\frac{1}{p}})\times[0,\,1)$. This leads to (i) by Theorem \ref{complete1}.

(ii) By a simple computation, we have
$$\Psi^{\ast}(x,\,\xi)\Psi(x,\,\xi)=V^{\ast}(x,\,\xi)\Lambda^{\ast}(x,\,\xi)\Lambda(x,\,\xi)V(x,\,\xi)$$
and
$$\Lambda^{\ast}(x,\,\xi)\Lambda(x,\,\xi)=diag\left(|\lambda_0(x,\,\xi)|^2,\,|\lambda_1(x,\,\xi)|^2,\,\cdots,\,
|\lambda_{p-1}(x,\,\xi)|^2\right)$$
for $(x,\,\xi)\in [1,\,a^{\frac{1}{p}})\times[0,\,1)$. It follows that
$$A^2I_p\leq \Psi^{\ast}(x,\,\xi)\Psi(x,\,\xi)\leq B^2I_p$$
for a.e. $(x,\,\xi)\in [1,\,a^{\frac{1}{p}})\times[0,\,1)$ if and only if $A\leq|\lambda_s(x,\,\xi)|\leq B$ for $s\in \bb N_p$ and a.e. $(x,\,\xi)\in [1,\,a^{\frac{1}{p}})\times[0,\,1)$. This leads to (ii) by Theorem \ref{frame1}.
\epr

\br By Lemma \ref{matrixexp}, Theorem \ref{Density1} covers all complete ${\cal MD}$-systems and ${\cal MD}$-frames for $L^2(\bb R_+)$. \er

\br\label{rem5} By Theorem \ref{frame1}, the lower frame bound $A$ and upper frame bound $B$ of a frame ${\cal MD}(\psi,\,a,\,b)$ for $L^2(\bb R_+)$ always satisfy that
\bea\label{equa13}a^{\frac{q-1}{p}}A\leq B.\eea
This shows that there exist significant differences between ${\cal MD}$-frames for $L^2(\bb R_+)$ and usual Gabor and wavelet frames for $L^2(\bb R)$. The following argument tells us that there exists a tight frame ${\cal MD}(\psi,\,a,\,b)$ for $L^2(\bb R_+)$ only if $a=b$. Indeed, suppose ${\cal MD}(\psi,\,a,\,b)$ is a tight frame for $L^2(\bb R_+)$ with frame bound $A$. Then $A=B$ in (\ref{equa13}), and thus $q=1$. Also observing that $p\leq q$ is necessary for the completeness of ${\cal MD}(\psi,\,a,\,b)$, we have $p=q=1$, equivalently, $a=b$.
\er

\br\label{rem6} Theorem \ref{Density1} shows that, for arbitrary positive numbers $A$ and $B$ satisfying $a^{\frac{q-1}{p}}A\leq B$, we can construct an ${\cal MD}$-frame ${\cal MD}(\psi,\,a,\,b)$ for $L^2(\bb R_+)$ with frame bounds $A$ and $B$. Indeed, take $\lambda_s$ with $s\in \bb N_p$ in Theorem \ref{Density1} such that $a^{\frac{q-1}{2p}}\sqrt{A}\leq\left|\lambda_s(x,\,\xi)\right|\leq\sqrt{B}$ for $s\in \bb N_p$ and a.e. $(x,\,\xi)\in [1,\,a^{\frac{1}{p}})\times[0,\,1)$. Then ${\cal MD}(\psi,\,a,\,b)$ is a frame for $L^2(\bb R_+)$ with frame bounds $A$ and $B$ by Theorem \ref{frame1}.
\er

Collecting Theorems \ref{complete1}-\ref{Density1} and Remark \ref{new}, we obtain the following density theorem:

\bt \label{Density} Let $a$ and $b$ be as in the general setup. Then the following are equivalent:

{\rm (i) }$\log_ba\leq 1$.

{\rm (ii) }There exists $\psi\in L^2(\bb R_+)$ such that ${\cal MD}(\psi,\,a,\,b)$ is complete in $L^2(\bb R_+)$.

{\rm (iii) }There exists $\psi\in L^2(\bb R_+)$ such that ${\cal MD}(\psi,\,a,\,b)$ is a frame for $L^2(\bb R_+)$.
\et

\bs
Finally, we conclude this paper by the following conjecture.\\
\textbf{Conjecture}.\label{rem7} In Theorem \ref{Density}, $\log_ba$ is required to be a rational number. This is a technical condition in all our arguments. We conjecture that, for general $a,\,b>1$, (i), (ii) and (iii) are equivalent to each other.

\end{document}